\documentclass[11pt]{amsart} \usepackage{latexsym} \usepackage{amsfonts}

\usepackage{amsmath,amsthm,amssymb}

\newtheorem{theor}{Theorem}

\newtheorem{thm}{Theorem}[section] 
\newtheorem{lem}[thm]{Lemma} 
\newtheorem{cor}[thm]{Corollary} 
 
\newtheorem{prop}[thm]{Proposition} 
\theoremstyle{definition} 
\newtheorem{defn}[thm]{Definition} 
 
\newtheorem{conv}[thm]{Convention} 
\newtheorem{rem}[thm]{Remark}

\newtheorem{quest}[thm]{Question} 
 
\def\strutdepth{\dp\strutbox} 
\def \ss{\strut\vadjust{\kern-\strutdepth \sss}} 
\def \sss{\vtop to \strutdepth{ 
\baselineskip\strutdepth\vss\llap{$\diamondsuit\;\;$}\null}} 
 
\def\strutdepth{\dp\strutbox} 
\def \sst{\strut\vadjust{\kern-\strutdepth \ssss}} 
\def \ssss{\vtop to \strutdepth{ 
\baselineskip\strutdepth\vss\llap{$\spadesuit\;\;$}\null}} 
 
\def\strutdepth{\dp\strutbox} 
\def \ssh{\strut\vadjust{\kern-\strutdepth \sssh}} 
\def \sssh{\vtop to \strutdepth{ 
\baselineskip\strutdepth\vss\llap{$\heartsuit\;\;$}\null}}

\newcommand{\R}{\mathbb R}

\newcommand{\Pro}{\mathbb P}

\begin{document} 
 
\author[Ilya Kapovich]{Ilya Kapovich}
 
\address{\tt Department of Mathematics, University of Illinois at
  Urbana-Champaign, 1409 West Green Street, Urbana, IL 61801, USA
  \newline http://www.math.uiuc.edu/\~{}kapovich/} \email{\tt
  kapovich@math.uiuc.edu}
 
\author[Martin Lustig]{Martin Lustig} \address{\tt Math\'ematiques
  (LATP), Universit\'e Paul C\'ezanne -Aix Marseille III, av. Escadrille
  Normandie-Ni\'emen, 13397 Marseille 20, France} \email{\tt
  Martin.Lustig@univ.u-3mrs.fr}
 
\thanks{The first author was supported by the NSF grant DMS\#0404991 and
  by the Humboldt Foundation Research Fellowship}
 
\title[Minimal sets and equivariant incompatibility]{The actions of
  $Out(F_k)$ on the boundary of Outer space and on the space of
  currents: minimal sets and equivariant incompatibility}

\begin{abstract} 
  We prove that for $k\ge 5$ there does not exist a continuous map
  $\partial CV(F_k)\to\mathbb PCurr(F_k)$ that is either
  $Out(F_k)$-equivariant or $Out(F_k)$-anti-equivariant. Here $\partial
  CV(F_k)$ is the ``length-function" boundary of Culler-Vogtmann's Outer
  space $CV(F_k)$, and $\mathbb PCurr(F_k)$ is the space of
  projectivized geodesic currents for $F_{k}$.

  We also prove that, if $k\ge 3$, for the action of $Out(F_k)$ on
  $\mathbb PCurr(F_{k})$ and for the diagonal action of $Out(F_k)$ on
  the product space $\partial CV(F_k)\times \mathbb PCurr(F_k)$ there
  exist unique non-empty minimal closed $Out(F_k)$-invariant sets.

  Our results imply that for $k\ge 3$ any continuous
  $Out(F_k)$-equivariant embedding of $CV(F_k)$ into $\mathbb
  PCurr(F_k)$ (such as the Patterson-Sullivan embedding) produces a new
  compactification of Outer space, different from the usual
  ``length-function" compactification $\overline{CV(F_k)}=CV(F_k)\cup
  \partial CV(F_k)$.
\end{abstract}

\subjclass[2000]{Primary 20F65, Secondary 30F60, 37B, 57M}
 
\keywords{geodesic currents, free groups, outer space}
 
\maketitle

\section{Introduction}\label{intro}

For a free group $F$ of finite rank $k\ge 2$
Culler-Vogtmann's~\emph{Outer space}~\cite{CV} $CV(F)$ is a fundamental
object for studying the group $Out(F)$ and the properties of individual
automorphisms of $F$. Outer space is a close cousin of the Teichm\"uller
space $\mathcal T(S_g)$ of a hyperbolic surface $S_g$. The analogy
between the action of the mapping class group on Teichm\"uller space and
the action of $Out(F)$ on $CV(F)$ is an important source of mathematical
and philosophical inspiration in the study of $Out(F)$. However, it is
well-understood that Outer space $CV(F)$ is more complicated than
Teichm\"uller space. The elements of $CV(F)$ are minimal free and
discrete actions of $F$ on $\mathbb R$-trees with quotient metric graphs
of volume one.  Any such $F$-action on a tree $T$ defines a translation
length function $|| \cdot ||_{T}: F \to \R_{\geq 0}$ and thus a point in
$\R^{F}$. This gives an embedding $CV(F) \subset \Pro\R^{F}$, and its
closure $\overline{CV(F)}$ is the ``length function'' compactification
of $CV(F)$, see \S 3 below. It is known that the elements of
$\overline{CV(F)}$ are precisely the projective classes of \emph{very
  small} minimal actions of $F$ on $\mathbb R$-trees~\cite{BF93, CL}.

The recent work~\cite{CHL3,KKS,Ka,Ka1,Ka2} shows that there is another
fundamental structure related to $Out(F)$, namely the space $Curr(F)$ of
\emph{geodesic currents} on $F$, that is, the space of positive Radon
$F$-invariant and flip-invariant measures on the space $\partial^2 F$ of
pairs of distinct points of $\partial F$.  Here the ``flip map" is
$\sigma: \partial^2 F\to \partial^2 F$, $\sigma: (\xi,\zeta)\mapsto
(\zeta,\xi)$. The corresponding projectivized space with respect to
scalar multiples is denoted $\mathbb PCurr(F)$. There is a natural
action of $Out(F)$ on $Curr(F)$ that quotients through to the action of
$Out(F)$ on $\mathbb PCurr(F)$.  The advantage of the projectivized
space $\mathbb PCurr(F)$ is that it is compact.

In the study of Teichm\"uller space geodesic currents are natural and
important objects, as elucidated, in particular, in the work of
Bonahon~\cite{Bo86,Bo88}. Thus the points of the Thurston boundary
$\partial \mathcal T(S_g)$ of Teichm\"uller space are geodesic
laminations on $S_g$ equipped with transverse measures. When pulled back
to the universal cover of $S_g$, these transverse measures lift to
$\pi_1(S_g)$-invariant measures on $\partial^{2} \mathbb
H^2=\partial^{2} \pi_1(S_g)$, that is, geodesic currents. These
transverse measures are important, in particular, because they are used
to define the metric on an $\mathbb R$-tree dual to a geodesic
lamination.

It turns out that for a free group $F$ geodesic currents and $\R$-trees
are naturally ``transversal objects". There is a canonical
$Out(F)$-equivariant \emph{intersection form} (see \cite{Ka,Ka1, Lu})
\[ 
I: cv(F)\times Curr(F)\to\mathbb R
\] 
where $cv(F)$ is the space of all (i.e. with arbitrary covolume)
$\R$-trees equipped with a free and discrete minimal action of $F$. This
intersection form has some important properties in common with Bonahon's
notion of an intersection number between geodesic currents on a
hyperbolic surface. In particular, for every $\R$-tree $T\in cv(F)$ and
every non-trivial $g\in F$ we have
\[ 
I(T, \mu_g)=|| g ||_{T},
\] 
where $\mu_g$ is the ``counting" or ``rational" current corresponding to
the conjugacy class $[g]$ (see Definition~\ref{defn:rational} below for
the precise definition of $\mu_g$), and $|| g ||_{T}$ denotes the
translation length on $T$ of the element $g \in F$.

The study of $Curr(F)$ has already led to some new results about the
geometry and dynamics of free group automorphisms. For example, it is
proved~\cite{Ka} that for every $\phi\in Aut(F)$ and every free basis
$A$ of $F$ the \emph{conjugacy distortion spectrum $I_A(\phi)$ of
  $\phi$} is a rational interval $I_A(\phi)=[\lambda_1,\lambda_2]\cap
\mathbb Q$ where by definition
\[ 
I_A(\phi):=\{\frac{||\phi(g)||_A}{||g||_A}: g\in F, g\ne 1\}.
\] 
Moreover, the extremal distortions $\lambda_1,\lambda_2$ actually belong
to $I_A(\phi)$, that is they are realized as distortions of some
non-trivial conjugacy classes. The geodesic currents approach also
provided a theoretical explanation for the experimental results
regarding the behavior of Whitehead's algorithm~\cite{Fr,Ka2,KSS}. The
work of Coulbois, Hilion and Lustig~\cite{CHL1,CHL2,CHL3} explored the
idea of an $\mathbb R$-tree dual to a measured geodesic lamination in
the context of free groups and some new and unexpected behavior was
uncovered.  Other results regarding currents and free group
automorphisms can be found in \cite{Fr,KKS,Ka,Ka1,Ka2,KLSS,KKS,KN,Ma}.

Since both $\overline{CV(F)}$ and $\mathbb PCurr(F)$ are fundamental and
intimately connected compact spaces endowed with natural
$Out(F)$-actions, it is interesting to understand how the dynamical
properties of these actions are related to each other.  There is an
additional reason to be interested in such questions, which is motivated
by what happens for hyperbolic surfaces.  Bonahon proved that if $S_g$
is a compact oriented hyperbolic surface and $\mathcal T(S_g)$ is the
Teichm\"uller space of $S_g$, then the Liouville map $L: \mathcal
T(S_g)\to\mathbb PCurr(\pi_1(S_g))$ is a topological embedding
equivariant with respect to the action of the mapping class group
$Mod(S_g)$ of $S_g$. Moreover, it turns out that the map $L$ extends to
a homeomorphism (which is necessarily $Mod(S_g)$-equivariant) from the
Thurston compactification $\overline{\mathcal T(S_g)}$ to the closure of
the image of $L$. It is natural to ask if there is an analogue of this
result for free groups, where the Thurston compactification of
Teichm\"uller space is replaced by the length function compactification
$\overline{CV(F)}$ of Outer space. It turns out that the answer to this
question is negative.  The reason for this is very general and is given
by the following theorem that we establish in this paper:

\begin{theor}\label{A} 
  Let $F$ be a free group of finite rank $k\ge 5$. Then there exists no
  continuous $Out(F)$-equivariant map $\partial CV(F)\to \mathbb
  PCurr(F)$. Similarly, there is no continuous $Out(F)$-anti-equivariant
  map $\partial CV(F)\to \mathbb PCurr(F)$.
\end{theor} 
 
There are several equivalent descriptions of the Liouville map
$L:\mathcal T(S_g)\to\mathbb PCurr(\pi_1(S_g))$ mentioned above.  One of
such descriptions involves characterizing the Liouville current as the
Patterson-Sullivan current corresponding to the hyperbolic structure. It
turns out that this characterization generalizes to the case of free
groups. In~\cite{KN} Kapovich and Nagnibeda proved that the
Patterson-Sullivan map $\mathcal P: CV(F)\to\mathbb PCurr(F)$ is an
$Out(F)$-equivariant topological embedding. Theorem~\ref{A} implies
that, provided $F$ has rank at least five, $\mathcal P$ does not extend
to a continuous map from $\overline{CV(F)}$.
 
\medskip
 
While proving Theorem~\ref{A} we obtain some new information about the
dynamics of the action of $Out(F)$ on $\mathbb PCurr(F)$.  In~\cite{Ma}
Reiner Martin introduces the following subset of $\mathbb PCurr(F)$: the
set $\mathcal M^{\mathbb PCurr}\subseteq \mathbb PCurr(F)$ is defined as
the closure in $\mathbb PCurr(F)$ of the set of projectivized rational
currents $[\mu_g]$ corresponding to all the primitive elements $g$ of
$F$. Although suggested by the name, Reiner Martin does not prove or
conjecture that $\mathcal M^{\mathbb PCurr}$ is the unique smallest
non-empty closed $Out(F)$-invariant subset of $\mathbb PCurr(F)$. We
shall prove here that this is indeed the case, provided $F$ has rank
bigger or equal to three:

\begin{theor}\label{B} 
  Let $F$ be a free group of finite rank $k\ge 3$. Then $\mathcal
  M^{\mathbb PCurr}\subseteq \mathbb PCurr(F)$ is the unique smallest
  non-empty closed $Out(F)$-invariant subset of $\mathbb PCurr(F)$.
\end{theor}

Since $\mathcal M^{\mathbb PCurr}$ is infinite dimensional~\cite{Ma},
when $F$ has rank $k\ge 3$, and $\partial CV(F)$ is finite
dimensional~\cite{Skora,Steiner,GL}, Theorem~\ref{B} implies:
 
\begin{cor} 
  Let $F$ be a free group of finite rank $k\ge 3$, there is no
  $Out(F)$-equivariant (or $Out(F)$-anti-equivariant) topological
  embedding from $\partial CV(F)$ to $\mathbb PCurr(F)$.
\end{cor} 
 
Thus we see that in this case the closure of the image of any
$Out(F)$-equivariant topological embedding $CV(F)\to\mathbb PCurr(F)$
gives a new compactification of Outer space $CV(F)$, different from the
length function compactification $\overline{CV(F)}$. In particular, this
applies to the Patterson-Sullivan embedding $\mathcal P: CV(F)\to\mathbb
PCurr(F)$. Together Theorem~\ref{A} and Theorem~\ref{B} imply:

\begin{cor} 
  Let $F$ be a finitely generated free group of rank $k\ge 2$ and let
  $\mathcal P :CV(F)\to\mathbb PCurr(F)$ be the Patterson-Sullivan
  embedding.
 
  If $k\ge 3$ then $\mathcal P$ does not extend to a homeomorphism from
  $\overline{CV(F)}$ to the closure of the image of $\mathcal P$, and
  $\mathcal P$ does not induce a topological embedding $\partial
  CV(F)\to \mathbb PCurr(F)$.
 
  Moreover, if $k\ge 5$, then $\mathcal P$ does not extend to a
  continuous map from $\overline{CV(F)}$ to $\mathbb PCurr(F)$.
\end{cor} 
 
All of the above results reflect the fact that the dynamics of the
$Out(F)$ action on $\mathbb PCurr(F)$ is rather different from the
dynamics of the $Out(F)$ action on $\overline{CV(F)}$ and on $\partial
CV(F)$. In particular, one can expect to find new information about the
dynamical properties of free group automorphisms by considering the
action of $Out(F)$ on the space of geodesic currents.
 
We will describe briefly the idea of the proof of Theorem~\ref{A}.
Suppose that $\tau: \partial CV(F)\to \mathbb PCurr(F)$ is an
$Out(F)$-equivariant continuous map. First we compare the dynamics of
the action of simple Dehn twists on $\partial CV(F)$ and $\mathbb
PCurr(F)$. In both cases this dynamics turns out to be of a ``parabolic"
nature. Using this fact we show that $\tau$ must take the length
function $[T_D]\in \partial CV(F)$, corresponding to a simple Dehn twist
$D$, to the rational current $[\mu_b]\in \mathbb PCurr(F)$ where $b\in
F$ is the ``twistor" of $D$ (see Section~\ref{sect:twists} for the
definitions related to simple Dehn twists).  Next we use the results of
Cohen-Lustig~\cite{CL} about the dynamics of Dehn multi-twists on
$\partial CV(F)$. We exhibit two specific points $[T_{D_1}],[T_{D_2}]\in
\partial CV(F)$, corresponding to simple Dehn twists $D_1$ and $D_2$
with the twistors $b_1$ and $b_2$, and a multi-twist $\phi$ of $F$ with
the following properties. On the trees side, we have $\lim_{n\to\infty}
\phi^n[T_{D_1}]=\lim_{n\to\infty} \phi^n[T_{D_2}]$. On the other hand,
on the currents side, $\phi, b_1,b_2$ are chosen so that
$\lim_{n\to\infty} \phi^n [\mu_{b_1}]\ne\lim_{n\to\infty} \phi^n
[\mu_{b_2}]$. Since, by the first step, we must have $\tau
([T_{D_1}])=[\mu_{b_1}]$ and $\tau ([T_{D_2}])=[\mu_{b_1}]$, we get a
contradiction with the assumption that $\tau$ is continuous and
$Out(F)$-equivariant.

The ideas involved in the proof of Theorem~\ref{A} also yield:
 
\begin{theor}\label{C} 
  Let be a free group of finite rank $k\ge 3$.
  
  Then there exists a unique minimal non-empty closed $Out(F)$-invariant
  subset $\mathcal M^{2}$ of $\partial CV(F)\times \mathbb PCurr(F)$.
  Moreover, this set $\mathcal M^{2}$ is equal to the closure of all
  points of $\partial CV(F)\times \mathbb PCurr(F)$ of the form $([T_D],
  [\mu_b])$ where $D$ is a simple Dehn twist of $F$ and $b$ is the
  twistor of $D$.
\end{theor} 
 
It turns out that in the context of Theorem~\ref{A} there are no
continuous equivariant or anti-equivariant maps going in the other
direction:

\begin{theor}\label{D} 
  Let $F$ be a free group of finite rank $k\ge 3$. Then there does not
  exist a continuous $Out(F)$-equivariant (or $Out(F)$-anti-equivariant)
  map $\mathcal M^{\mathbb PCurr}\to\partial CV(F)$. Hence there does
  not exist a continuous $Out(F)$-equivariant (or
  $Out(F)$-anti-equivariant) map $\mathbb PCurr(F)\to\partial CV(F)$.
\end{theor} 
 
The proof of Theorem~\ref{D} turns out to be somewhat easier than that
of Theorem~\ref{A}. First, we show, again by exploring the parabolic
dynamics, that if $\tau:\mathcal M^{\mathbb PCurr}\to\partial CV(F)$ is
a map as in Theorem~\ref{D} then for every simple Dehn twist $D$ of $F$
with twistor $b$ we have $\tau([\mu_b])=[T_D]$. Then, since $k\ge 3$, it
is easy two produce two simple Dehn twists $D$ and $D'$ with the same
twistor $b$ such that $[T_D]\ne [T_{D'}]$, yielding a contradiction.

\section{Basic definitions} 
 
We will give only a quick review of the basic concepts related to 
geodesic currents on free groups. We refer the reader to 
\cite{Ka,Ka1} for a detailed treatment of this topic. 
 
\begin{conv} 
  For the remainder of the paper, unless specified otherwise, let $F$ be
  a finitely generated free group of rank $k\ge 2$. We will denote by
  $\partial F$ the hyperbolic boundary of $F$ in the sense of the theory
  of word-hyperbolic groups. Since $F$ is free, $\partial F$ can also be
  viewed as the space of ends of $F$ with the standard ends-space
  topology.
 
  Thus $\partial F$ is a topological space homeomorphic to the Cantor
  set. We will also denote
\[ 
\partial^2 F:=\{(\zeta,\xi): \zeta,\xi\in \partial F \text{ and }
\zeta\ne \xi\}.
\] 
 
Denote by $\sigma:\partial^2 F\to \partial^2 F$ the \emph{flip} map
$\sigma:(\zeta,\xi)\mapsto (\xi,\zeta)$ for $(\zeta,\xi)\in \partial^2
F$.
 
\end{conv} 
 
\begin{defn}[Geodesic Currents] 
  Let $F$ be a free group of finite rank $k\ge 2$. A \emph{geodesic
    current} on $F$ is a positive Radon measure on $\partial^2 F$ that
  is $F$-invariant and $\sigma$-invariant. We denote the space of all
  geodesic currents on $F$ by $Curr(F)$.
 
  The space $Curr(F)$ comes equipped with a weak topology that in this
  case can be characterized as follows: for $\nu_n,\nu\in Curr(F)$ we
  have $\displaystyle\lim_{n\to\infty}\nu_n=\nu$ if and only if for
  every two disjoint closed-open sets $S,S'\subseteq \partial F$ we have
  $\displaystyle\lim_{n\to\infty}\nu_n(S\times S')=\nu(S\times S')$.
\end{defn} 
 
Note that the above definition and notations are consistent with
\cite{CHL1,CHL2,CHL3} but are a little different from those used in
\cite{Ka1}. Namely, in \cite{Ka1} $Curr(F)$ denotes all $F$-invariant
positive Radon measure on $\partial^2 F$. The subspace of those such
measures that are also $\sigma$-invariant is denoted in \cite{Ka1} by
$Curr_s(F)$.

\begin{defn}[Projectivized Geodesic Currents] 
  For two non-zero geodesic currents $\nu_1, \nu_2\in Curr(F)$ we say
  that $\nu_1$ is {\em projectively equivalent} to $\nu_2$, denoted
  $\nu_1\sim \nu_2$, if there exists a non-zero scalar $r\in \mathbb R$
  such that $\nu_2=r \nu_1$. We denote \[\mathbb PCurr(F):=\{\nu\in
  Curr(F): \nu\ne 0\}/\sim\] and call it the \emph{space of
    projectivized geodesic currents on $F$}. Elements of $\mathbb
  PCurr(F)$ (that is, scalar equivalence classes of elements of
  $Curr(F)$) are called \emph{projectivized geodesic currents}. The
  space $\mathbb PCurr(F)$ is endowed with the quotient topology. We
  will denote the $\sim$-equivalence class of a non-zero geodesic
  current $\nu$ by $[\nu]$.
\end{defn} 
 
\begin{defn}[$Out(F)$-action] 
  Let $\phi\in Aut(F)$ be an automorphism. It is well known that $\phi$
  extends to a homeomorphism of $\partial F$ and hence induced a
  homeomorphism $\hat\phi$ of $\partial^2 F$.
 
  For any $\nu\in Curr(F)$ define a measure $\phi\nu$ on $\partial^2F$
  as a pull-back:
  $$(\phi\nu)(S):=\nu\big((\hat\phi)^{-1}(S)\big)$$
  for $S\subseteq
  \partial^2 F$. It can be shown~\cite{Ka1} that $\phi\nu\in Curr(F)$ is
  a geodesic current on $F$ and, moreover, $\phi\nu$ only depends on
  $\nu$ and on the outer automorphism class $[\phi]\in Out(F)$.  The map
  $Aut(F)\times Curr(F)\to Curr(F)$, $(\phi,\nu)\mapsto \phi\nu$ defines
  a continuous left action of $Aut(F)$ on $Curr(F)$ by linear
  transformations that factors through to a left continuous action of
  $Out(F)$ on $Curr(F)$. Moreover, this action commutes with scalar
  multiplication and hence defines a continuous left action of $Out(F)$
  on $\mathbb PCurr(F)$.
\end{defn}

\begin{defn}[Coordinates on $Curr(F)$] 
  Let $A=\{a_1,\dots, a_k\}$ be a free basis of $F$ and let $Cay(F,A)$
  be the Cayley graph of $F$ with respect to $A$. Thus $Cay(F,A)$ is a
  $2k$-regular tree.
 
  Let $\gamma$ be a non-empty geodesic segment in $Cay(F,A)$ that begins
  and ends at vertices of $Cay(F,A)$. The path $\gamma$ comes equipped
  with a \emph{label} $v$ which is a freely reduced word over $A$.
 
  We denote by $Cyl_A(\gamma)$ the set of all $(\xi,\zeta)\in \partial^2
  F$ such that the directed geodesic from $\xi$ to$\zeta$ in $Cay(F,A)$
  contains the segment $\gamma$. Thus $Cyl_A(\gamma)$ is a closed-open
  subset of $\partial^2 F$.
 
  Let $\nu\in Curr(F)$ be arbitrary and let $v\in F$ be a non-trivial
  freely reduced word. Since $\nu$ is $F$-invariant, if $\gamma$ is a
  segment in $Cay(F,A)$ with label $v$, then $\nu(Cyl_A(\gamma))$
  depends only on $\nu$ and $v$ but not on the choice of a lift $\gamma$
  of $v$.  Moreover, since $\nu$ is $\sigma$-invariant, we have
  $\nu(Cyl_A(\gamma))=\nu(Cyl_A(\gamma^{-1}))$.

  We denote
\[ 
(v; \nu)_A:=\frac{1}{2}(\nu(Cyl_A(\gamma))+\nu(Cyl_A(\gamma^{-1})))
=\nu(Cyl_A(\gamma)),
\] 
where $\gamma$ is any lift of $v$ to $Cay(F,A)$.
\end{defn} 
 
\begin{defn}[Rational currents]\label{defn:rational} 
  Let $g\in F$ be a non-trivial element that is not a proper power.
  Then there exist two distinct points $g_+,g_-\in \partial F$ such that
  $\lim_{n\to\infty} g^n=g_+$ and $\lim_{n\to\infty} g^{-n}=g_{-}$. We
  define
\[ 
\mu_g:= \sum_{h\in [g]}
\left(\delta_{(h_-,h_+)}+\delta_{(h_+,h_-)}\right),
\] 
where $[g]$ is the conjugacy class of $g$ in $F$, and $\delta_{(\zeta,
  \xi)}$ is the Dirac measure defined by the point $(\zeta, \xi) \in
\partial^{2}F$.
 
If $g\in F$ is a non-trivial element that is a proper power, we can
uniquely write $g=f^m$ where $m>1$ and $f$ is not a proper power. We set
$\mu_g:=m\mu_f$.
 
Thus for every non-trivial $g\in F$ we have $\mu_g\in Curr(F)$ is 
a geodesic current on $F$. We say that $\nu\in Curr(F)$ is a 
\emph{rational current} if $\nu$ has the form $\nu= s \mu_g$ for 
some $s>0$ and $g\in F$. 
\end{defn}

\begin{prop}\cite{Ma,Ka,Ka1} 
  The set of rational currents is dense in $Curr(F)$.
\end{prop}

\begin{rem}[Rational currents and cyclic words]\label{rem:cyclic} 
  Let $A$ be a free basis of $F$. It is often convenient to represent
  conjugacy classes in $F$ by ``cyclic words". A \emph{cyclic word} $w$
  over $A$ is a non-trivial cyclically reduced word in $A$ written on a
  circle clockwise without a specified base-point. If $v$ is a freely
  reduced word, a vertex on a cyclic word $w$ is an \emph{occurrence} of
  $v$ in $w$ if we can read $v$ in $w$ from this vertex going forward
  clockwise and without leaving the labelled circle. The number of
  occurrences of $v^{\pm 1}$ in $w$ is denoted $(v; w)_A$. Thus by
  definition $(v; w)_A=(v^{-1}; w)_A$. We denote the number of vertices
  in a cyclic word $w$ in the basis $A$ by $||w||_{A}$; this coincides
  with the word length of $w$ in $A^{\pm 1}$, since $w$ is (cyclically)
  reduced.
 
  Clearly there is a bijective correspondence between cyclic words over
  $A$ and non-trivial conjugacy classes in $F$.
 
  A simple but important observation~\cite{Ka,Ka1} says that if a cyclic
  word $w$ represents a conjugacy class $[g]$ in $F$ then for every
  freely reduced word $v$ over $A$ we have
\[ 
(v;\mu_g)_A:=(v; w)_A.
\] 
\end{rem}

\begin{defn}[Length]\label{defn:length} 
Let $A$ be a free basis of $F$ and let $\nu\in Curr(F)$ be a 
current. Denote 
\[ 
||\nu||_A:=\sum_{a\in A} (a; \nu)_A.
\] 
We call $||\nu||_A$ the \emph{length} of $\nu$ with respect to $A$.
\end{defn} 
 
Remark~\ref{rem:cyclic} implies that for any $a \in A$ the bracket $(a;
\mu_g)_A$ equals to the number of occurrences of $a$ or $a^{-1}$ in the
cyclically reduced word in $A^{\pm 1}$ that represents $g$. For every
non-trivial $g\in F$ we have $||\mu_g||_A=||g||_A$, the cyclically
reduced length of $g$ with respect to $A$. Note also that for any
current $\mu$ the length $||\mu||_A$ is precisely equal to the
intersection form $I(T_{A}, \mu)$, where $T_{A}$ denotes the Cayley tree
associated to the basis $A$, with all edges of length $1$.
 
We will need the following basic facts~\cite{Ka1}, where for any element
$g \in F$ we denote by $|g|_A$ the length of the reduced word in the
basis $A$ that represents $g$.
 
\begin{lem}\label{lem:basic} Let $A$ be a free basis of $F$ and let $\nu\in 
  Curr(F)$. Then for every freely reduced word $v\in F$ we have
 
\begin{enumerate} 
\item $(v; \nu)_A\ge 0$;
 
\item $(v; \nu)_A=\sum_{a\in A^{\pm 1}: |va|_A= |v|_{A}+1} (va; \nu)_A$;
 
\item $(v; \nu)_A=\sum_{a\in A^{\pm 1}: |av|_A=|v|_{A}+1} (av; \nu)_A$;
 
\item For every $m\ge 1$
\[ 
||\nu||_A=\frac{1}{2}\sum_{u\in F : |u|_{A}=m} (u; \nu)_A.
\] 

\item Let $\nu_n,\nu\in Curr(F)$. Then $\lim_{n\to\infty} \nu_n=\nu$ in
  $Curr(F)$ if and only if for every $v\in F$, $v\ne 1$ we have
\[
\lim_{n\to\infty} (v; \nu_n)_A=(v; \nu)_A.
\]

\item Let $\nu_n,\nu\in Curr(F)$ be nonzero currents. Then we have
  $$\lim_{n\to\infty} [\nu_n]=[\nu]\quad \text{ in } \quad \mathbb
  PCurr(F)$$
  if and only if
  $$\lim_{n\to\infty}\frac{\nu_n}{||\nu_n||_A}=\frac{\nu}{||\nu||_A}
  \quad \text{ in }\quad Curr(F).$$

\end{enumerate} 
 
\end{lem} 
 
\begin{cor}\label{cor:limit}
  Let $A$ be a free basis of $F$ and let $a\in A$. For nonzero currents
  $\nu_n\in Curr(F)$ we have
\[
\lim_{n\to\infty}  [\nu_n]=[\mu_a]
\]
if and only if
\[
\lim_{n\to\infty} \frac{(a;\nu_n)_A}{||\nu_n||_A}=1.
\]

\end{cor}

\begin{proof}
  The ``only if" direction is obvious. Suppose that $\lim_{n\to\infty}
  \frac{(a;\nu_n)_A}{||\nu_n||_A}=1$.  We claim that
  $\lim_{n\to\infty}\frac{\nu_n}{||\nu_n||_A}=\mu_a$ in $Curr(F)$.  Note
  that
\[
\big|\big| \frac{\nu_n}{||\nu_n||_A} \big|\big|_A=1.
\]

Definition~\ref{defn:length} implies that $\lim_{n\to\infty}
\frac{(x;\nu_n)_A}{||\nu_n||_A}=0$ for every $x\in A, x\ne a$. If $x\in
A$, $x\ne a$ then by parts (2) and (3) of Lemma~\ref{lem:basic} for
every freely reduced word $u$ involving $x^{\pm 1}$ we have $(u;
\nu_n)_A\le (x;\nu_n)_A$. Therefore for such $u$ we have
$\lim_{n\to\infty} \frac{(u;\nu_n)_A}{||\nu_n||_A}=0$. Hence, by part
(4) of Lemma~\ref{lem:basic} for every $m\ge 1$ we have
$\lim_{n\to\infty} \frac{(a^m;\nu_n)_A}{||\nu_n||_A}=1$. Thus we have
verified that for every nontrivial $v\in F$ we have $\lim_{n\to\infty}
\frac{(v;\nu_n)_A}{||\nu_n||_A}=(v; \mu_a)$. Hence
$\lim_{n\to\infty}\frac{\nu_n}{||\nu_n||_A}=\mu_a$ in $Curr(F)$ and
$\lim_{n\to\infty} [\nu_n]=[\mu_a]$ in $\mathbb PCurr(F)$, as required.

\end{proof}

Corollary~\ref{cor:limit} and Remark~\ref{rem:cyclic} easily imply:

\begin{cor}\label{cor:conv}
  Let $A$ be a free basis of $F$ and let $a\in A$. Then for every $g\in
  F$ we have

\[
\lim_{n\to\infty} [\mu_{a^ng}]=\lim_{n\to\infty} [\mu_{a^{-n}g}]=[\mu_a]
\]
in $\mathbb PCurr(F)$.
\end{cor}

\section{Conventions regarding 
  Outer space}
 
We refer the reader to~\cite{BF00,CL,CV,Gui,Vog} for detailed
information regarding Outer space and its boundary. We recall here
briefly some basic notions and facts.
 
Let $F$ be a free group of finite rank $k\ge 2$. We denote by $cv(F)$
the \emph{non-projectivized Outer space} of $F$, which consists of all
$\R$-trees $T$ equipped with an action of $F$ by isometries which is
minimal, free and discrete. Such an action is characterized (up to
$F$-equivariant isometries) by the associated translation length
function $|| \cdot ||_{T}: F \to \R$, and for the purpose of this paper
it is convenient to identify the point $T \in cv(F)$ with the
corresponding length function $|| \cdot ||_{T} \in \R^{F}$, i.e. $cv(F)
\subset \R^{F}$. Note that for every $T \in cv(F)$ the quotient by the
$F$-action is a finite connected {metric marked graph} without valence 1
vertices: this gives an alternative characterization of the points in
$cv(F)$.  A typical example is the Cayley tree $T_{A}$ associated to a
basis $A$ of $F$, which has as quotient a rose with $\#A$ leaves,
usually assumed to be of unit length.

We denote by $CV(F)$ the subset of $cv(F)$ corresponding to those
actions where the quotient metric graphs have volume one.
Alternatively, we can also view $CV(F)$ as the projectivized quotient
$\mathbb Pcv(F):=cv(F)/\sim$, where $\sim$ corresponds to scalar
multiplication. Namely, for every $T \in cv(F)$ there exists a unique
rescaled tree $\lambda T$ with the quotient of volume one. The space
$CV(F)$ is called the \emph{Outer space} of $F$, and the projective
class of a tree $T$ is denoted by $[T]$.
 
There is a \emph{left} action $Aut(F)\times cv(F)\to cv(F)$ given by $||
w ||_{\phi T} = ||\phi^{-1}(w) ||_{T}$, for all $w \in F$ and all $T \in
cv(F)$. This action leaves $CV(F)$ invariant.  Moreover, the subgroup
$Inn(F)$ of inner automorphisms of $F$ acts trivially on $cv(F)$ and the
actions of $Aut(F)$ on $cv(F)$ and $CV(F)$ factor through to left
actions of $Out(F)$ on $cv(F)$ and $CV(F)$ accordingly.

The space $cv(F) \subset \R^{F}$ is endowed with the weak topology of
pointwise convergence on finite subsets of $F$. The action of $Out(F)$
on $cv(F)$ is an action by homeomorphisms and $CV(F)\subseteq cv(F)$ is
a closed $Out(F)$-invariant subset.
 
Outer space $CV(F)=\mathbb Pcv(F) \subset \mathbb P \R^{F}$ can be
compactified with respect to weak convergence of projective classes of
length functions, and the result is denoted $\overline{CV(F)}$ and is
called the \emph{length function compactification} of $CV(F)$. The
difference $\overline{CV(F)}-CV(F)$ is denoted $\partial CV(F)$ and is
called the \emph{length-function boundary} or sometimes the {\em
  Thurston boundary} of $CV(F)$. It turns out~\cite{BF93, CL} that
$\overline{CV(F)}$ consists precisely of all \emph{very small} minimal
actions of $F$ on $\mathbb R$-trees.

\section{Minimal sets} 
 
Throughout this section assume that the rank $k$ of $F$ satisfies $k\ge
3$. As Guirardel showed in~\cite{Gui}, there exists a unique minimal
non-empty closed $Out(F)$-invariant subset $\mathcal M^{CV}\subseteq
\partial CV(F)$. This set obviously has the property that for every
$x\in \mathcal M^{CV}$ the orbit $Out(F)x$ is dense in $\mathcal
M^{CV}$. Moreover, Guirardel describes in detail some of the points that
belong to $\mathcal M^{CV}$. We shall only need the following basic
fact~\cite{Gui}:
 
\begin{prop}\label{prop:gui} 
  There exists a point $[T]\in \partial CV(F)$ corresponding to a free
  action of $F$ on an $\mathbb R$-tree $T$ such that $[T]\in \mathcal
  M^{CV}$. Moreover, since the $Out(F)$-orbit of $[T]$ is dense in
  $\mathcal M^{CV}$, the set of points of $\mathcal M^{CV}$
  corresponding to free actions of $F$ is dense in $\mathcal M^{CV}$.
\end{prop} 
 
Specifically, one can choose $[T]$ in Proposition~\ref{prop:gui} to be
the attracting fixed point in $\partial CV(F)$ of any outer automorphism
$\alpha$ of $F$ which is irreducible with irreducible powers and does
not have periodic conjugacy classes. Levitt and Lustig~\cite{LL} proved
that such automorphisms have ``North-South" dynamics on $\partial CV(F)$
(and indeed on the entire space $\overline{CV(F)}$) and hence, clearly,
their attracting fixed points must belong to $\mathcal M^{CV}$. It is
also known that these attracting fixed points correspond to free
actions.
 
\medskip

Similarly but not quite analogously, Reiner Martin~\cite{Ma} introduces
the following notion in the context of currents:
 
\begin{defn}[Minimal set in $\mathbb PCurr(F)$] 
  Set $\mathcal M^{\mathbb PCurr}\subseteq \mathbb PCurr(F)$ to be the
  closure in $\mathbb PCurr(F)$ of the set
\[ 
\{[\mu_g] : g\in F \text{ is a primitive element} \}.
\] 
\end{defn} 
 
Here by a primitive element we mean an element that belongs to some free
basis of $F$.  Recall that, if $\phi\in Aut(F)$ and $g\in F$ is
non-trivial, then $\phi [\mu_g]=[\mu_{\phi(g)}]$. Note that for any
primitive element $g$ the set $\mathcal M^{\mathbb PCurr}$ is the
closure of the $Out(F)$-orbit of $[\mu_g]$. Thus $\mathcal M^{\mathbb
  PCurr}$ is a closed $Out(F)$-invariant subset of $\mathbb PCurr(F)$.
We will need the following basic result of Reiner Martin~\cite{Ma}. For
completeness we include a proof.
 
\begin{prop}\label{prop:rmartin} 
  Let $F$ be a free group of finite rank $k\ge 3$. Then the space
  $\mathcal M^{\mathbb PCurr}$ is infinite dimensional.
\end{prop}

\begin{proof} 
 
We will denote $F=F(a,b,c,\dots )$ and $F':=F(a,b)$. There is a 
natural topological embedding (see~\cite{Ka1} for proof) \[ \iota: 
\mathbb PCurr(F')\to \mathbb PCurr(F)\] with the property that for 
every $g\in F'$ we have $\iota([\mu_g'])=\mu_g$ where $\mu_g'\in 
Curr(F')$ and $\mu_g\in Curr(F)$ are the rational currents 
corresponding to $g$ considered as the element of $F'$ and $F$ 
accordingly.

We claim that $\mathcal M^{\mathbb PCurr}$ contains a copy
$\iota(\mathbb PCurr(F'))$ of $\mathbb PCurr(F')$. Since rational
currents are dense in $\mathbb PCurr(F')$, to establish the claim it
suffices to show that for every $g\in F'$ we have $\mu_g\in \mathcal
M^{\mathbb PCurr}$. Let $u=u(a,b)\in F'$ be an arbitrary non-trivial
freely reduced word. Then for every $n\ge 1$ the element $g_n:=cu^n$ is
primitive in $F$. Thus $[\mu_{g_n}]\in \mathcal M^{\mathbb PCurr}$. It
is easy to see that $\lim_{n\to\infty} [\mu_{g_n}]=[\mu_u]\in \mathbb
PCurr(F)$. Thus indeed $\mathcal M^{\mathbb PCurr}$ contains the image
of the topological embedding $\iota$. Since $\mathbb PCurr(F')$ is
infinite dimensional~\cite{Ma,Ka,Ka1} it now follows that so is
$\mathcal M^{\mathbb PCurr}$.
 
\end{proof} 
 
\section{Simple Dehn twists and their action on $\mathbb PCurr(F)$}\label{sect:twists} 
 
\begin{defn}\label{defn:D} 
  Let $F$ be a free group of finite rank $k\ge 2$ with a free basis
  $A=\{a_1,\dots, a_k\}$. Denote $a=a_1, b=a_2$. Let $D:F\to F$ be the
  automorphism defined as $D(a)=ab$, $D(a_i)=a_i$ for $2\le i\le k$. We
  will call $D$ a \emph{simple Dehn twist} with respect to $A$ and say
  that $b$ is the \emph{twistor} of $D$.
\end{defn} 
 
\begin{defn}[Splitting corresponding to a Dehn twist]\label{defn:T_D} 
  Let $D$ be a simple Dehn twist of $F$ as in Definition~\ref{defn:D}.
 
  Denote $b_1:=aba^{-1}\in F$. Note that the elements $b_1,b,a_3,\dots,
  a_k$ freely generate a subgroup of $F$ that we denote by $F_0$.
 
  Consider the following HNN-extension decomposition of $F$:
\[ 
F=\langle F_0, a| a b a^{-1}=b_1 \rangle.
\] 
 
Let $T_D\in cv(F)$ be the Bass-Serre tree corresponding to this
HNN-extension. We call $T_D$ \emph{the tree corresponding to the simple
  Dehn twist $D$}.
\end{defn}

\begin{defn}\label{defn:crit}[Critical set] 
  Let $D$ be the simple Dehn twist with respect to a free basis $A$ with
  twistor $b$, as in Definition~\ref{defn:D}. Denote
\[ 
Y:=\{[\nu]: 0\ne \nu\in Curr(F),
\frac{1}{2}(a;\nu)_A=(ab;\nu)_A=(ab^{-1};\nu)_A\}\subseteq \mathbb
PCurr(F).
\] 
 
We call $Y$ the \emph{critical set} of $D$ in $\mathbb PCurr(F)$.
\end{defn}

\begin{conv}\label{conv:D} 
  Until the end of this section, unless specified otherwise, we assume
  that $D$ is a simple Dehn twist of $F$ with twistor $b$ as in
  Definition~\ref{defn:D}.
\end{conv} 
 
We will need the following fact~\cite{CL}:
 
\begin{prop}\label{prop:D} 
  Let $[T]\in \overline {CV(F)}$ be such that $||b||_T > 0$. Then
\[ 
\lim_{n\to\infty} D^n [T]=\lim_{n\to\infty} D^{-n} [T]=[T_D].
\] 
\end{prop}

Recall that in our coordinate notation we have
$(a_i;\nu)_A=(a_i^{-1};\nu)_A$ for $i=1,\dots, k$ and for every $\nu\in
Curr(F)$.
 
\begin{lem}\label{D-basic} Let $\nu\in Curr(F)$ be a current. We have: 
\begin{enumerate} 
\item $(a_i; D\nu)_A=(a_i; \nu)_A$ for $i\ne 2$.

\item $(b; D\nu)_A=(b; \nu)_A+(a;\nu)_A-2(ab^{-1}; \nu)_A$.

\item $(ab^{-1}; D\nu)_A=(ab^{-2};\nu)_A+(ab^{-1}a^{-1}; \nu)_A$ and
  hence $(ab^{-1}; D\nu)_A\le (ab^{-1},\nu)_A$.
 
\item $||D\nu||_A-||\nu||_A=(b; D\nu)_A-(b;\nu)_A=(a;\nu)_A-2(ab^{-1};
  \nu)_A$. Consequently, $||D\nu||_A>||\nu||_A$ if and only if
  $(a;\nu)_A>2(ab^{-1}; \nu)_A$.
\end{enumerate} 
 
\end{lem} 
 
\begin{proof} 
  Parts (1), (2) and (3) are straightforward for cyclic words and hence
  for rational currents. Therefore they follow for arbitrary currents by
  continuity since rational currents are dense in $Curr(F)$. Clearly,
  part (4) follows from parts (1) and (2).
 
\end{proof} 
 
\begin{lem}\label{D-dynamics} 
  Let $\nu\in Curr(F)$ and suppose that $||D\nu||_A>||\nu||_A$. Then
  $||D^2\nu||_A>||D\nu||_A$ and $||D^2\nu||_A-||D\nu||_A\ge
  ||D\nu||_A-||\nu||_A$.
\end{lem} 
 
\begin{proof} 
 
  By assumption and by Lemma~\ref{D-basic} (4) we have
  $||D\nu||_A-||\nu||_A=(a;\nu)_A-2(ab^{-1}; \nu)_A>0$.

  By Lemma~\ref{D-basic} (1) and (3) we have $(a;D\nu)_A=(a;\nu)_A$ and
  $(ab^{-1};D\nu)_A\le (ab^{-1};\nu)_A$.
 
  Hence
\begin{gather*} 
  ||D^2\nu||_A-||D\nu||_A=(a;D\nu)_A-2(ab^{-1}; D\nu)_A\ge\\
  (a;\nu)_A-2(ab^{-1}; \nu)_A=||D\nu||_A-||\nu||_A>0
\end{gather*} 
and the statement of the lemma follows.

\end{proof} 
 
\begin{cor}\label{cor:D} 
  Let $\nu\in Curr(F)$ and suppose that $(a;\nu)_A>2(ab^{-1}; \nu)_A$.
  Then we have:
 
\begin{gather*} 
  \lim_{n\to\infty} ||D^n\nu||_A=\infty, \quad
  \lim_{n\to\infty}\frac{1}{||D^n\nu||_A} D^n \nu =\mu_b\quad \text{
    and}\\
  \lim_{n\to\infty} [D^n\nu]=[\mu_b] \quad \text{ in }\quad\mathbb
  PCurr(F).
\end{gather*} 
\end{cor} 
 
\begin{proof} 
  Lemma~\ref{D-basic} (4) and Lemma~\ref{D-dynamics} imply that
  $||D\nu||_A>||\nu||_A$ and that for $n\ge 1$ we have \[||D^n
  \nu||_A\ge ||\nu||_A+ n(||D\nu||_A-||\nu||_A).\] Hence
  $\lim_{n\to\infty} ||D^n\nu||=\infty$. For every $a_i\ne b$ we have
  $(a_i,D^n\nu)_A=(a_i, \nu)_A$. Hence for $a_i\ne b$ one obtains
\[ 
\lim_{n\to\infty} \frac{(a_i,D^n\nu)_A}{||D^n\nu||_A}=0.
\] 
 
From the definition of $||D^n\nu||_A$ it follows
\[ 
\lim_{n\to\infty} \frac{(b,D^n\nu)_A}{||D^n\nu||_A}=1.
\] 
Therefore by Corollary~\ref{cor:limit} one has $\lim_{n\to\infty}
[D^n\nu]=[\mu_b]$ in $\mathbb PCurr(F)$, as claimed.
 
\end{proof}

\begin{cor}\label{cor:D1} 
  Let $\nu\in Curr(F)$ be a non-zero current such that $[\nu]\not\in Y$.
  Then in $\mathbb PCurr(F)$ one has either $\lim_{n\to\infty}
  [D^n\nu]=[\mu_b]$ or $\lim_{n\to\infty} [D^{-n}\nu]=[\mu_b]$.
\end{cor} 
 
\begin{proof} 
  We claim that at least one of the following holds: $(a;\nu)_A>2(ab;
  \nu)_A$ or $(a;\nu)_A>2(ab^{-1}; \nu)_A$. Suppose not. Then
  $\frac{1}{2}(a;\nu)_A\le (ab; \nu)_A$ and $\frac{1}{2}(a;\nu)_A\le
  (ab^{-1}; \nu)_A$. Since $[\nu]\not\in Y$, at least one of these
  inequalities must be strict and hence $(a;\nu)_A<(ab; \nu)_A+(ab^{-1};
  \nu)_A$. On the other hand, Lemma~\ref{lem:basic} (4) implies that
  $(a;\nu)_A\ge (ab; \nu)_A+(ab^{-1}; \nu)_A$, yielding a contradiction.
 
  Thus we have indeed either $(a;\nu)_A>2(ab; \nu)_A$ or
  $(a;\nu)_A>2(ab^{-1}; \nu)_A$. If the former holds, then
  Corollary~\ref{cor:D} applied to $D^{-1}$ gives $\lim_{n\to\infty}
  [D^{-n}\nu]=[\mu_b]$. If the latter inequality holds then
  Corollary~\ref{cor:D} applied to $D$ we have $\lim_{n\to\infty}
  [D^{n}\nu]=[\mu_b]$.
 
\end{proof}

We can now prove that $\mathcal M^{\mathbb PCurr}$ is indeed the
smallest non-empty closed $Out(F)$-invariant subset of $F$ assuming that
$F$ has rank at least three.

\begin{thm}[Minimality of $\mathcal M^{\mathbb PCurr}$]\label{thm:min} 
  Let $F$ be a free group of rank $k\ge 3$. Let $S\subseteq \mathbb
  PCurr(F)$ be a non-empty closed $Out(F)$-invariant subset. Then one
  has:
  $$\mathcal M^{\mathbb PCurr}\subseteq S.$$
\end{thm} 
 
\begin{proof} 
  Fix a free basis $A=\{a_1,a_2,a_3,\dots, a_k\}$ of $A$ and denote
  $a=a_1$, $b=a_2$, $c=a_3$.

  Since $S$ is non-empty, there exists a non-zero current $\nu\in
  Curr(F)$ such that $[\nu]\in S$.

  There exists $i$ such that $(a_i; \nu)_A > 0$ and without loss of
  generality we may assume that $i=1$, so that $(a; \nu)_A > 0$.

  Let $D$ be the simple Dehn twist with twistor $b$ defined as $D(a)=ab$
  and $D(a_i)=a_i$ for $i\ge 2$. Let $Y\subseteq \mathbb PCurr(F)$ be
  the critical set of $D$ as in Definition \ref{defn:crit}.
 
  Suppose first that $[\nu]\not\in Y$. Then by Corollary~\ref{cor:D1}
  either $\lim_{n\to\infty} [D^n\nu]=[\mu_b]$ or $\lim_{n\to\infty}
  [D^{-n}\nu]=[\mu_b]$. In either case, since $S$ is closed and
  $Out(F)$-invariant, it follows that $[\mu_b]\in S$ and hence $\mathcal
  M^{\mathbb PCurr}\subseteq S$.

  Suppose now that $[\nu]\in Y$, that is
\[ 
\frac{1}{2}(a;\nu)_A= (ab^{-1}; \nu)_A=(ab; \nu)_A
\] 
By part (4) of Lemma~\ref{lem:basic} this implies that
$(ac;\nu)_A=(ac^{-1};\nu)_A=0$.

Consider now the simple Dehn twist $D'$ of $F$ defined as $D'(a)=ac$,
$D'(a_i)=a_i$ for $2\le i\le k$. Thus $c$ is the twistor of $D'$. Since
$(ac^{-1};\nu)_A=0$ and $(a;\nu)_A>0$, we have
$(a;\nu)_A-2(ac^{-1};\nu)_A>0$. Thus by Corollary~\ref{cor:D} applied to
$D'$ we have $\lim_{n\to\infty} [(D')^n\nu]=[\mu_c]$. Again, since $S$
is closed and $Out(F)$-invariant it follows that $[\mu_c]\in S$ and
hence $\mathcal M^{\mathbb PCurr}\subseteq S$.
 
\end{proof}

\begin{cor}\label{cor:noemb} 
  Let $F$ be a free group of finite rank $k\ge 3$. Then there does not
  exist an $Out(F)$-equivariant topological embedding of $\partial
  CV(F)$ to $\mathbb PCurr(F)$.
 
  Similarly, there exists no $Out(F)$-anti-equivariant topological
  embedding of $\partial CV(F)$ to $\mathbb PCurr(F)$.
\end{cor} 
 
\begin{proof} 
  Suppose $j:\partial CV(F)\to \mathbb PCurr(F)$ is an
  $Out(F)$-equivariant (or $Out(F)$-anti-equivariant) topological
  embedding. Since $\partial CV(F)$ is compact, this implies that the
  image of $j$ is a closed $Out(F)$-invariant subset of $\mathbb
  PCurr(F)$. Hence by Theorem~\ref{thm:min} the image of $j$ contains
  the set $\mathcal M^{\mathbb PCurr}$. By
  Proposition~\ref{prop:rmartin} the set $\mathcal M^{\mathbb PCurr}$ is
  infinite dimensional. On the other hand, $\partial CV(F)$ is finite
  dimensional~\cite{Skora,Steiner,GL}, yielding a contradiction with the
  fact that $j$ is a homeomorphism from $\partial CV(F)$ onto its image
  $j(\partial CV(F))$.
 
\end{proof} 
 
In particular, Corollary~\ref{cor:noemb} implies that if $\tau: CV(F)\to
\mathbb PCurr(F)$ is an $Out(F)$-equivariant topological embedding then
$\tau$ does not extend to a homeomorphism from the length-function
compactification $\overline{CV(F)}=CV(F)\cup \partial CV(F)$ to the
closure of the image of $\tau$.

\section{Rigid points} 
 
\begin{conv} 
  For the remainder of this section we assume that $F$ is a free group
  of rank $k\ge 3$ and that $\tau: \partial CV(F)\to \mathbb PCurr_S(F)$
  is a continuous map that is either $Out(F)$-equivariant or
  $Out(F)$-anti-equivariant.
\end{conv} 
 
\begin{prop}\label{prop:rat} 
  For any simple Dehn twist $D\in Aut(F)$ with twistor $b \in F$ we have
  $\tau([T_D])=[\mu_b]$.
\end{prop}

\begin{proof} 
  Let $A$ be the basis of $F$ used to define $D$ as in
  Definition~\ref{defn:D}.  Let $Y\subseteq \mathbb PCurr(F)$ be the
  critical set of $D$, as defined in Definition~\ref{defn:crit}. Thus
\[ 
Y=\{[\nu]: \nu\in Curr(F), \nu\ne 0,
\frac{1}{2}(a;\nu)_A=(ab;\nu)_A=(ab^{-1};\nu)_A\}.
\]

Note that $Y$ is a closed non-empty subset of $\mathbb PCurr(F)$ and
that $Y$ contains the set
\[ 
Y_0:=\{[\nu]\in \mathbb PCurr(F): \nu\in Curr(F), \nu\ne 0,
(a;\nu)_A=0\}.
\] 
Moreover, $[\mu_a]\not\in Y$ and hence $\mathcal M^{\mathbb
  PCurr}\not\subseteq Y$.
 
We claim that there exists $[T]\in \mathcal M^{CV}$ such that $||b||_T>
0$ and $\tau([T])\not\in Y$. Indeed, suppose not, so that for every
$[T]\in \mathcal M^{CV}$ with $||b||_T> 0$ we have $\tau([T])\in Y$. By
Proposition~\ref{prop:gui} the set of length functions corresponding to
free actions is dense in $\mathcal M^{CV}$ and for each of them $b$ has
non-zero translation length.  Hence $\tau(\mathcal M^{CV})\subseteq Y$
since $Y$ is closed. The space $\mathcal M^{CV}$ is compact and
$Out(F)$-invariant. Hence $\tau(\mathcal M^{CV})$ is a non-empty compact
(and thus closed) subset of $\mathbb PCurr(F)$ that is
$Out(F)$-invariant. Therefore by Theorem~\ref{thm:min} $\mathcal
M^{\mathbb PCurr}\subseteq \tau(\mathcal M^{CV})\subseteq Y$, yielding a
contradiction with the earlier conclusion that $\mathcal M^{\mathbb
  PCurr}\not\subseteq Y$.
 
Thus indeed there exists $x=[T]\in \mathcal M^{CV}$ such that $||b||_T>
0$ and $\tau(x)\not\in Y$.

Therefore, by Corollary~\ref{cor:D1} we have
 
\[ 
\text{ either } \lim_{n\to\infty} D^n \tau(x)=[\mu_b] \text{ or }
\lim_{n\to\infty} D^{-n} \tau(x)=[\mu_b].
\] 
Suppose the former holds (the other case is symmetric) and that we have
$\displaystyle \lim_{n\to\infty} D^n \tau(x)=[\mu_b]$.
 
If $\tau$ is $Out(F)$-equivariant, then by Proposition~\ref{prop:D} we
have $[T_D]=\lim_{n\to\infty} D^n x$, and hence the continuity of $\tau$
implies
\[\tau([T_D])=\lim_{n\to\infty} \tau(D^n x)=\lim_{n\to\infty} 
D^n\tau(x)=[\mu_b].\]
 
If $\tau$ is $Out(F)$-anti-equivariant, then, again by
Proposition~\ref{prop:D} we have $[T_D]=\lim_{n\to\infty} D^{-n} x$ and
hence \[\tau([T_D])=\lim_{n\to\infty} \tau(D^{-n} x)=\lim_{n\to\infty}
D^n\tau(x)=[\mu_b].\]
 
Thus $\tau([T_D])=[\mu_b]$, as required.
 
\end{proof}

\section{Non-existence of equivariant maps} 
 
\begin{conv} 
  Throughout this section let $F$ be a free group of finite rank $k\ge
  5$ with a free basis $A=\{a_1,\dots, a_k\}$. We denote $a=a_1$,
  $b=a_2$, $c=a_3$, $d=a_4$ and $e=a_5$.
 
  Let $\phi$ denote the automorphism of $F$ defined as $\phi(a)=ab$,
  $\phi(e)=ed$ and $\phi(a_i)=a_i$ for $i\ne 1,5$.
\end{conv} 
 
\begin{defn}[The length function $\Delta(\rho,\theta)$.] 
  Denote $b_1=aba^{-1}$, $d_1=ede^{-1}$ and put $H:=\langle b,b_1,c,
  d,d_1, a_6,\dots, a_k\rangle\le F$. (Note that $H$ is freely generated
  by these elements).
 
  Consider the following HNN-extension splitting of $F$ with stable
  letters $a,e$ and the base $H$:
 
\[ 
F=\langle H,a,e| aba^{-1}=b_1, ede^{-1}=d_1\rangle.
\] 
 
Let $\mathbb A$ be the graph of groups corresponding to this splitting.
Let $\rho>0,\, \theta>0$ be positive numbers.
 
We denote by $\Delta(\rho,\theta)\in cv(F)$ the hyperbolic length
function on $F$ coming from $\mathbb A$ where in $\mathbb A$ the
loop-edge labelled $a$ is given length $\rho$ and the loop-edge labelled
$e$ is given length $\theta$.
 
\end{defn} 
 
The following result is due to Cohen-Lustig (see \cite{CL}, Theorem
13.2).
 
\begin{prop}\label{prop:ch} 
  Let $[T]\in \partial CV(F)$ be such that $||b||_T> 0$ and $||d||_T>
  0$.
 
  Then
\[ 
\lim_{n\to\infty} \phi^n([T])=\lim_{n\to\infty}
\phi^{-n}([T])=[\Delta(||b||_T, ||d||_T)]
\] 
in $\partial CV(F)$.
\end{prop}

\begin{thm}\label{thm:main} 
  Let $F$ be a free group of rank $k\ge 5$. Then there does not exist a
  continuous $Out(F)$-equivariant map $\tau:\partial CV(F)\to \mathbb
  PCurr(F)$. Similarly, there exists no continuous
  $Out(F)$-anti-equivariant map $\partial CV(F)\to \mathbb PCurr(F)$.
\end{thm} 
\begin{proof} 
 
  Suppose $\tau:\partial CV(F)\to \mathbb PCurr(F)$ is a continuous map
  that is $Out(F)$-equivariant.

  Consider the free bases $A'=\{a',b',c',d',e',a_6,\dots, a_k\}$ and
  $A''=\{a'',b'',$ $c'',d'',e'',a_6,\dots, a_k\}$ of $F$ defined via
\[ 
a=a', e=e', b=b'a'c', d=a'b'd', c=a'b'
\] 
and\[ a=a'', e=e'', b=b''a''c'', d=a''b''d'', c=e''b''.
\] 
Note that we have $b'=a^{-1}c$ and $b''=e^{-1}c$.
 
Let $D'$ and $D''$ be automorphisms of $F$ defined as follows. We have
$D'(a')=(a'b')$ and $D'$ fixes all other elements of $A'$.  Similarly
$D''(a'')=a''b''$ and $D''$ fixes all other elements of $A''$. Let
$T_{D'}$ and $T_{D''}$ be the length functions on $F$ defined as in
Definition~\ref{defn:T_D}.
 
Since $b=b'a'c'$ and $d=a'b'd'$, we have $||b||_{T_{D'}}=1$ and
$||d||_{T_{D'}}=1$. Similarly, $b=b''a''c''$ and $d=a''b''d''$, so that
$||b||_{T_{D''}}=1$ and $||d||_{T_{D''}}=1$. Therefore by
Proposition~\ref{prop:ch} we have
\[ 
\lim_{n\to\infty} \phi^n [T_{D'}]=\lim_{n\to\infty} \phi^n
[T_{D''}]=[\Delta(1,1)]\in \partial CV(F).
\] 
 
On the other hand, by Proposition~\ref{prop:rat} we have
$\tau([T_{D'}])=[\mu_{b'}]=[\mu_{a^{-1}c}]$ and
$\tau([T_{D''}])=[\mu_{b''}]=[\mu_{e^{-1}c}]$.
 
By definition of $\phi$ for $n\ge 1$ we have
$\phi^n(a^{-1}c)=b^{-n}a^{-1}c$ and hence
$\phi^n[\mu_{b'}]=\phi^n[\mu_{a^{-1}c}]=[\mu_{b^{-n}a^{-1}c}]$.
Similarly, for $n\ge 1$ we have $\phi^n(e^{-1}c)=d^{-n}e^{-1}c$ and
hence $\phi^n[\mu_{b''}]=\phi^n[\mu_{e^{-1}c}]=[\mu_{d^{-n}e^{-1}c}]$.
Therefore, by Corollary~\ref{cor:conv}, we see that
\[ 
\lim_{n\to\infty}\phi^n[\mu_{b'}]=[\mu_b]\in \mathbb PCurr(F)
\] 
and
\[ 
\lim_{n\to\infty}\phi^n[\mu_{b''}]=[\mu_d]\in \mathbb PCurr(F).
\] 
 
Since $[\mu_b]\ne [\mu_d]$, we have
\begin{gather*} 
  \lim_{n\to\infty} \phi^n [T_{D'}]=\lim_{n\to\infty} \phi^n
  [T_{D''}] \text{  in  } \partial CV(F) \text{  but }\\
  \lim_{n\to\infty} \phi^n \tau[T_{D'}]\ne \lim_{n\to\infty} \phi^n
  \tau[T_{D''}] \text{ in } \mathbb PCurr(F),
\end{gather*} 
yielding a contradiction with the assumption that $\tau$ is continuous
and $Out(F)$-equivariant.
 
The same proof shows that there does not exist a continuous
$Out(F)$-anti-equivariant map $\tau_0: \partial CV(F)\to \mathbb
PCurr(F)$. The only changes needed to be made in the argument are the
following. Using Proposition~\ref{prop:rat} we still conclude that
$\tau_0([T_{D'}])=[\mu_{b'}]=[\mu_{a^{-1}c}]$ and
$\tau_0([T_{D''}])=[\mu_{b''}]=[\mu_{e^{-1}c}]$. As before, we have
\[ 
\lim_{n\to\infty} \phi^n [T_{D'}]=\lim_{n\to\infty} \phi^n [T_{D''}]
\text{ in } \partial CV(F).
\] 
By anti-equivariance, $\tau_0(\phi^n x)=\phi^{-n} \tau_0(x)$ for every
$x\in \partial CV(F)$. For $n\ge 1$ we have
$\phi^{-n}(a^{-1}c)=b^na^{-1}c$, $\phi^{-n}(e^{-1}c)= d^ne^{-1}c$ and so
$\phi^{-n}([\mu_{a^{-1}c}])=[\mu_{b^na^{-1}c}]$ and
$\phi^{-n}([\mu_{e^{-1}c}])=[\mu_{d^ne^{-1}c}]$. Therefore, again by
Corollary~\ref{cor:conv}, we have
\[ 
\lim_{n\to\infty}\phi^{-n}([\mu_{a^{-1}c}])=[\mu_b] \text{ and }
\lim_{n\to\infty}\phi^{-n}([\mu_{e^{-1}c}])=[\mu_d], \] and, since
$[\mu_b]\ne [\mu_d]$ we get a contradiction with the continuity of
$\tau_0$ as before.
 
\end{proof} 
 
\begin{rem} 
  The proof of Theorem~\ref{thm:main} actually implies that there exists
  neither an $Out(F)$-equivariant nor an $Out(F)$-anti-equivariant
  continuous map from $\mathcal M^{CV}$ to $\mathbb PCurr(F)$. Indeed,
  the same holds with $\mathcal M^{CV}$ replaced by any non-empty closed
  $Out(F)$-invariant (or anti-invariant) subspace of $\partial CV(F)$.
\end{rem}

\section{The minimal set in $\partial CV(F)\times \mathbb 
  PCurr(F)$}
 
\begin{defn} 
  Let $F$ be a free group of finite rank $k\ge 3$. Let $\mathcal M^{2}$
  be the closure in $\partial CV(F)\times \mathbb PCurr(F)$ of all
  points of the form $([T_D], [\mu_g])$, where $D$ is a simple Dehn
  twist of $F$ with respect to any base, and $g \in F$ is the twistor of
  $D$.  Recall that $g$ is always a primitive element, according to the
  conventions in this paper, see Definition \ref{defn:D}.
\end{defn}

Note that by definition $\mathcal M^{2}$ is an $Out(F)$-invariant subset
of $\partial CV(F)\times \mathbb PCurr(F)$. In fact, if $D$ is any
simple Dehn twist of $F$ with twistor $b$, then $\mathcal M^{2}$ is the
closure of the $Out(F)$-orbit of the point $([T_D], [\mu_b])$.

\begin{thm}\label{thm:min-prod} 
  Let be a free group of finite rank $k\ge 3$. Then $\mathcal M^{2}$ is
  the unique minimal non-empty closed $Out(F)$-invariant subset of
  $\partial CV(F)\times \mathbb PCurr(F)$.
\end{thm} 
\begin{proof} 
  Let $S\subseteq \partial CV(F)\times \mathbb PCurr(F)$ be a non-empty
  closed $Out(F)$-invariant subset of $\partial CV(F)\times \mathbb
  PCurr(F)$. We need to show that $\mathcal M^{2}\subseteq S$.
 
  It is clear that the coordinate projections of $\mathcal M^{2}$ to the
  spaces $\partial CV(F)$ and $\mathbb PCurr(F)$ are compact (and hence
  closed) $Out(F)$-invariant sets and therefore they contain $\mathcal
  M^{CV}$ and $\mathcal M^{\mathbb PCurr}$ accordingly.

  Choose a simple Dehn twist $D$ of $F$ with a twistor $b$. Let
  $Y\subseteq \mathbb PCurr(F)$ be the critical set of $D$, defined in
  Definition~\ref{defn:crit}.
 
  We claim that there exists a point $([T],[\nu])\in S$ such that
  $||b||_T> 0$ and $\nu\not\in Y$. Suppose not. Then for every point
  $([T],[\nu])\in S$ either $||T||_b=0$ or $\nu\in Y$. Choose
  $([T],[\nu])\in S$ such that $T$ corresponds to a free action of $F$.
  (Such a point exists since the first coordinate projection of $S$
  contains $\mathcal M^{CV}$). Then $\phi [T]$ corresponds to a free
  action for every $\phi\in Out(F)$. For every $\phi\in Out(F)$ we have
  $\phi([T],[\nu])=(\phi[T], \phi[\nu])\in S$ and thus $\phi\nu\in Y$.
  Therefore the closure $C$ of the orbit $Out(F) ([T],[\nu])$ is a
  subset of $\partial CV(F)\times Y$. The second coordinate projection
  of $C$ is a closed $Out(F)$-invariant subset of $\mathbb PCurr(F)$
  that is contained in $Y$. This is a contradiction to
  Theorem~\ref{thm:min}, since $Y$ does not contain $\mathcal M^{\mathbb
    PCurr}$.
 
  Thus the claim is verified and there exists a point $([T],[\nu])\in S$
  such that $||b||_T> 0$ and $\nu\not\in Y$. By definition of $Y$ one
  has either $\frac{1}{2}(a;\nu)_A > (ab^{-1};\nu)_A$ or
  $\frac{1}{2}(a;\nu)_A > (ab;\nu)_A$. Suppose the former (the other
  case is symmetric). Then by Proposition~\ref{prop:D} and
  Corollary~\ref{cor:D} we have
\[ 
\lim_{n\to\infty} D^n ([T], [\nu])=([T_D],[\mu_b]).
\] 
Since $S$ is closed and $Out(F)$-equivariant, it follows that
$([T_D],[\mu_b])\in S$. Hence the closure of the $Out(F)$-orbit of
$([T_D],[\mu_b])$, that is the set $\mathcal M^{2}$, is also contained
in $S$, as claimed.
 
\end{proof} 
 
\section{Maps from currents to the boundary of 
  Outer space}
 
\begin{lem}\label{lem:dense} 
  Let $F$ be a free group of finite rank $k\ge 3$ and let $D$ be a
  simple Dehn twists of $F$ with twistor $b$. Let $Y\subseteq \mathbb
  PCurr(F)$ be the critical set of $D$. Then $\mathcal M^{\mathbb
    PCurr}-Y$ is dense in $\mathcal M^{\mathbb PCurr}$.
\end{lem} 
 
\begin{proof} 
  Recall that rational currents corresponding to primitive elements are
  dense in $\mathcal M^{\mathbb PCurr}$, by the definition of the
  latter.  Hence it suffices to show that if $g$ is a primitive element
  with $[\mu_g]\in Y$ then $[\mu_g]$ can be approximated by elements of
  $\mathcal M^{\mathbb PCurr}-Y$.
 
  Let $A=\{a,b,c,\dots \}$ be the free basis of $F$ that appears in the
  definition of the simple Dehn twist $D$. Suppose $g$ is a primitive
  element such that $[\mu_g]\in Y$. By definition of $Y$ it follows
  that, when expressed as a reduced word in $A^{\pm 1}$, the element $g$
  involves an even number of occurrences of $a^{\pm 1}$. After replacing
  $g$ by its conjugate if necessary, we may assume that $g$ is
  represented by a cyclically reduced word in $A$.

  Choose a free basis $B$ of $F$ containing $g$. By looking at the
  abelianization of $F$ we see that there must exist an element $f\in
  B$, $f\ne g^{\pm 1}$ such that $f$, when expressed in $A$, involves an
  odd number of occurrences of $a^{\pm 1}$. Thus $(a; \mu_f)_A$ is odd
  and hence by definition of $Y$ we have $[\mu_f]\not\in Y$. Moreover,
  for every integer $n\ge 1$ we have that $(a; \mu_{g^nf})_A$ is odd and
  hence $[\mu_{g^nf}]\not\in Y$.  On the other hand $g^nf$ is primitive
  for every integer $n$ and
\[ 
\lim_{n\to\infty} [\mu_{g^nf}]=[\mu_g].
\] 
Since $[\mu_{g^nf}]\in \mathcal M^{\mathbb PCurr}-Y$, the lemma is
proved.
 
\end{proof}

\begin{prop}\label{prop:back} 
Let $F$ be a free group of finite rank $k\ge 3$. Let $\tau: 
\mathcal M^{\mathbb PCurr}\to\partial CV(F)$ be a continuous map 
that is either $Out(F)$-equivariant or $Out(F)$-anti-equivariant. 
 
Then for every simple Dehn twist $D$ of $F$ with twistor $b$ we have
$\tau([\mu_b])=[T_D]$.
 
\end{prop} 
 
\begin{proof} 
  Let $Y\subseteq \mathbb PCurr(F)$ be the critical set of $D$.
 
  We claim that there exists a point $[\mu]\in \mathcal M^{\mathbb
    PCurr}-Y$ such that for $[T]=\tau([\mu])$ we have $||b||_T> 0$.
  Suppose not. Then, since by Lemma~\ref{lem:dense} the set $\mathcal
  M^{\mathbb PCurr}-Y$ is dense in $\mathcal M^{\mathbb PCurr}$ and
  since $\tau$ is continuous, it follows that the image of the map
  $\tau$ is contained in the set $\{[T]\in \partial CV(F): ||b||_T=0\}$.
  Since $\tau$ is either $Out(F)$-equivariant or
  $Out(F)$-anti-equivariant, it follows that the set $\{[T]\in \partial
  CV(F): ||b||_T=0\}$ contains a closed $Out(F)$-invariant subset of
  $\partial CV(F)$ and hence it contains the set $\mathcal M^{CV}$. This
  contradicts the fact that $\mathcal M^{CV}$ has elements corresponding
  to free actions.
 
  Thus indeed, there exists a point $[\mu]\in \mathcal M^{\mathbb
    PCurr}-Y$ such that for $[T]=\tau([\mu])$ we have $||b||_T> 0$.
  Then $\lim_{n\to \infty} D^n([T])=\lim_{n\to \infty}
  D^{-n}([T])=[T_D]$ and either $\lim_{n\to\infty} D^n([\mu])=[\mu_b]$
  or $\lim_{n\to\infty} D^{-n}([\mu])=[\mu_b]$.  By assumptions on
  $\tau$ it follows that $\tau([\mu_b])=[T_D]$, as claimed.
 
\end{proof} 
 
\begin{thm}\label{thm:back} 
  Let $F$ be a free group of finite rank $k\ge 3$. Then there does not
  exist a continuous map $\tau: \mathcal M^{\mathbb PCurr}\to\partial
  CV(F)$ that is either $Out(F)$-equivariant or
  $Out(F)$-anti-equivariant.
\end{thm} 
 
\begin{proof} 
 
  Suppose $\tau: \mathcal M^{\mathbb PCurr}\to\partial CV(F)$ is a
  continuous map that is either $Out(F)$-equivariant or
  $Out(F)$-anti-equivariant.
 
  Let $A=\{a_1,\dots, a_k\}$ be a free basis of $F$. Denote $a_1=a$,
  $a_2=b$ and $a_3=c$. Let $D$ be the simple Dehn twist with twistor $b$
  defined as $D(a)=ab$ and $D(a_i)=a_i$ for $i\ge 2$.
 
  Consider the free basis $A'=\{a_1',\dots, a_k'\}$ where $a_3'=ca$ and
  $a_i'=a_i$ for $i\ne 3$. Let $D'$ be the simple Dehn twist defined as
  $D'(a_1')=a_1'a_2'=ab$ and $D'(a_i')=a_i'$ for $i\ge 2$. Then $D'$ has
  twistor $a_2'=b$. Since both $D$ and $D'$ are simple Dehn twists with
  the same twistor $b$, Proposition~\ref{prop:back} implies that
  $\tau([\mu_b])=[T_D]=[T_{D'}]$.  However, for $c':=a_3'=ca$ we have
  $||c'||_{T_D}=1$ and $||c'||_{T_{D'}}=0$. Hence $[T_D]\ne [T_{D'}]$,
  yielding a contradiction.
 
\end{proof}

Since $\mathcal M^{\mathbb PCurr}$ is a closed $Out(F)$-invariant subset
of $\mathbb PCurr(F)$, Theorem~\ref{thm:back} immediately implies:

\begin{cor} 
  Let $F$ be a free group of finite rank $k\ge 3$. Then there does not
  exist a continuous map $\tau: \mathbb PCurr(F)\to\partial CV(F)$ that
  is either $Out(F)$-equivariant or $Out(F)$-anti-equivariant.
\end{cor} 

\section{Outlook}

The following question seems natural, and a positive answer would be
useful for several applications, some indicated below. In this section
we assume that $F$ is a finitely generated free group of rank $k\ge 3$.

\begin{quest}\label{prob:ext}
  Let $z(F)$ be the closure of the (non-projectivized) Outer space
  $cv(F)$ in the ambient vector space $\R^{F}$, given by the canonical
  embedding $cv(F) \subset \R^{F}$ via the translation length function
  as explained in \S 3.
  
  Does the intersection form $I:cv(F)\times Curr(F)\to\mathbb R$ admit a
  continuous extension $\overline{I}: z(F)\times Curr(F)\to\mathbb R$ ?
\end{quest}

Note that if such $\overline I$ exists, it must be $Out(F)$-invariant
since the intersection form $I$ is $Out(F)$-invariant. Similarly,
$\overline I$ will have to be $\mathbb R$-homogeneous with respect to
the first argument and linear with respect to the second argument.

The second author has produced a preliminary preprint~\cite{Lu} which
addresses this question and sketches the proof of a positive answer.
Working out some of the details still needs to be done.

If a continuous extension $\overline I$ does exist, then the notion of
having zero intersection number makes sense for the elements of
$\overline {CV(F)}$ and $\mathbb PCurr(F)$. Namely, for $[T]\in
\overline {CV(F)}$ and $[\nu]\in \mathbb PCurr(F)$ we say that
$([T],[\nu])\in I_0$ if $\overline I(T,\nu)=0$. It is easy to see that
$I_0$ has to be closed and $Out(F)$-invariant.

\begin{quest}\label{prob:prod}
  Let $F$ be free of finite rank $k\ge 3$. Recall that $\mathcal M^{2}$
  is the smallest non-empty closed $Out(F)$-invariant subset of
  $\partial CV(F)\times \mathbb PCurr(F)$. Is it true that
\[
\mathcal M^{2}=I_0\cap \left(\mathcal M^{CV}\times \mathcal M^{\mathbb
    PCurr}\right) \, ?
\]
\end{quest}

We do not know whether the answer to the inclusion "$\supset$" should
expected to be positive. For the inclusion "$\subset"$ we will now give
an argument, based on the assumption that Question \ref{prob:ext} has a
positive answer:

We already know that for a simple Dehn twist $D$ of $F$ corresponding to
$A$ with twistor $b\in A$ we have $ ([T_D],[\mu_b])\in \mathcal M^{2}$.
Let $T\in cv(F)$ be the action of $F$ on its Cayley graph with respect
to $A$.  It is known and easy to see directly that $\frac{1}{n} D^n
T\to_{n\to\infty} T_D$ in $z(F)$. Also, we have $D^n \mu_a=\mu_{D^n
  a}=\mu_{ab^n}$, $||D^n\mu_a||_A=||ab^n||_A=n+1$ and
$\lim_{n\to\infty}\frac{1}{n+1} D^n\mu_a =\mu_b$.  Therefore
\begin{gather*}
  I(\frac{1}{n} D^n T, \frac{1}{n+1}D^n \mu_a)=\frac{1}{n (n+1)} I(D^n
  T,D^n \mu_a
  )=\\
  \frac{1}{n (n+1)} I(T,\mu_a)\to_{n\to\infty} 0.
\end{gather*}
Hence the presumed continuity of $\overline I$ implies that
$([T_D],[\mu_b])\in I_0$. By equivariance of $\overline I$ we get $\phi
([T_D],[\mu_b])\in I_0$ for every $\phi\in Out(F)$. Since the
$Out(F)$-orbit of $([T_D],[\mu_b])$ is dense in $\mathcal M^{2}$, it
follows that $\mathcal M^{2}\subseteq I_0$ and therefore
\[
\mathcal M^{2}\subseteq I_0\cap \left(\mathcal M^{CV}\times \mathcal
  M^{\mathbb PCurr}\right).
\]

In this context we would also point to the work in progress of the
second author with Coulbois and Hilion \cite{CHL4} where it is shown
that the minimal set $\mathcal M^{2}$ is contained in a subset $\mathcal
L^{2}$ of $\mathcal M^{CV}\times \mathcal M^{\mathbb PCurr}$ that is
described in terms of algebraic laminations: A pair $(T, \mu)$ defines a
point of $\mathcal L^{2}$ if the {\em dual algebraic lamination}
$L^{2}(T)$ associated to $T$ in \cite{CHL2} contains the {\em support}
$L^{2}(\mu)$ associated to $\mu$ in \cite{CHL3}.  Recent results lead
the second author to the belief that $\mathcal L^{2}$ agrees with
$I_0\cap \left(\mathcal M^{CV}\times \mathcal M^{\mathbb PCurr}\right)$.

\begin{quest}
  Kapovich and Nagnibeda~\cite{KN} proved that the Patterson-Sullivan
  map $\mathcal P: CV(F)\to\mathbb PCurr(F)$ is a continuous
  $Out(F)$-equivariant topological embedding. Hence the closure
  $\overline{image(\mathcal P)}$ of the image of $\mathcal P$ is a
  closed $Out(F)$-invariant set. Hence by Theorem~\ref{thm:min}
  $\overline{image(\mathcal P)}$ contains the minimal set $\mathcal
  M^{\mathbb PCurr}$.
  
  Is it true that
\[
\overline{image(\mathcal P)}= image(\mathcal P)\cup \mathcal M^{\mathbb
  PCurr} \, ?
\]
\end{quest}

We would like to finish this section with a question of more speculative
character:

\begin{quest}
  For any point $[T] \in CV(F)$, is the set $\mathcal P^{2}$ of
  accumulation points of the $Out(F)$-orbit of the pair $([T], [\mathcal
  P(T)])$ strictly smaller than $I_0\cap \left(\mathcal M^{CV}\times
    \mathcal M^{\mathbb PCurr}\right)$ (or than $\mathcal L^{2}$) ?  Is
  it perhaps equal to $\mathcal M^{2}$ ?
\end{quest}

Note that for any current $0 \ne \mu \in Curr(F)$ and for any $[T] \in
CV(F)$ the set of accumulation points of the $Out(F)$-orbit of the pair
$([T], [\mu])$ is contained in $\partial CV \times \mathcal PCurr(F)$,
and, since it is closed, $Out(F)$-invariant and non-empty, it must
contain $\mathcal M^{2}$. We do not know whether it depends on the
choice of $T$ and $\mu$, and we do also not know whether it can ever be
strictly bigger than $\mathcal M^{2}$.

\end{document}